\documentclass{article}
\usepackage{graphicx,xcolor} 
\usepackage{amsmath,amssymb,amsthm}
\usepackage{setspace}
\usepackage[backend=biber,style=alphabetic]{biblatex}
\usepackage[nameinlink,noabbrev]{cleveref}
\usepackage{aliascnt}
\title{Bounding and Approximating Kostka Numbers}
\author{Jonathan Leake and Zirou Zhang}

\begin{document}
\theoremstyle{definition}

\newtheorem{definition}{Definition}[section]

\newaliascnt{theorem}{definition}
\newtheorem{theorem}[theorem]{Theorem}
\aliascntresetthe{theorem}

\newaliascnt{lemma}{definition}
\newtheorem{lemma}[lemma]{Lemma}
\aliascntresetthe{lemma}

\newaliascnt{remark}{definition}
\newtheorem{remark}[remark]{Remark}
\aliascntresetthe{remark}

\newaliascnt{example}{definition}
\newtheorem{example}[example]{Example}
\aliascntresetthe{example}

\newaliascnt{algorithm}{definition}
\newtheorem{algorithm}[algorithm]{Algorithm}
\aliascntresetthe{algorithm}

\newcommand{\MainTheoremBody}{%
Let $\lambda$ be a partition and let $\mu\in\mathbb{N}^n$ with $n\ge l(\lambda)$. If $K_{\lambda,\mu}>0$, then
\begin{align*}
\operatorname{Cap}_{\mu}(s_{\lambda}(x_1,\ldots,x_n))
\ge K_{\lambda,\mu}
&\ge \left(\frac{n}{e(n+|\lambda|)}\right)^n
\operatorname{Cap}_{\mu}(s_{\lambda}(x_1,\ldots,x_n)) 
\end{align*}
Note that $\left(\frac{n}{e(n+|\lambda|)}\right)^n \geq e^{-n-|\lambda|}$, giving an approximation factor of $e^{-n-|\lambda|}$.
Moreover, there is an optimization problem $(P_\mu)$ with optimal value $g_\mu^*$ such that
\[
\operatorname{Cap}_{\mu}(s_{\lambda}(x_1,\ldots,x_n))
= e^{g_\mu^*}.
\]
Given $\epsilon>0$, there is an algorithm that returns a $g_\mu^\circ$ satisfying
\[
g_\mu^* \le g_\mu^\circ \le g_\mu^*+\epsilon,
\]
with running time polynomial in $n$, $|\lambda|$, and $\log(1/\epsilon)$.
}
\maketitle

\begin{abstract}

    Kostka numbers, the coefficients of Schur polynomials, count semistandard Young tableaux of prescribed shape and weight. We give a deterministic algorithm for approximating Kostka numbers within a simply exponential factor in the number of variables and the size of the partition. Our approach combines Gurvits capacity bounds for log-concave polynomials with an algorithm for approximating the Gurvits capacity of a Schur polynomial. For a Schur polynomial in \(n\) variables indexed by \(\lambda\), we show that the logarithm of the capacity associated with any nonzero coefficient can be approximated to additive error \(\epsilon\) in time polynomial in \(n\), \(|\lambda|\), and \(\log(1/\epsilon)\). The analysis establishes an explicit bounded optimization region and efficient approximate first-order oracles, and handles boundary weights through a factorization into smaller Schur polynomials.
\end{abstract}

\section{Introduction}\label{sec:introduction}

\textit{Schur polynomials} are an important class of symmetric functions indexed by partitions. Their coefficients, known as the \textit{Kostka numbers}, have deep connections to combinatorics, representation theory, and algebraic geometry. In particular, Kostka numbers count semistandard Young tableaux of a fixed shape and weight. They also appear as weight multiplicities in irreducible polynomial representations of general linear groups. However, the direct computation of nonzero Kostka numbers is known to be \#$P$-complete \cite{narayanan2005computationkostkanumberslittlewoodrichardson}.

The \textit{capacity method}, developed by Gurvits, uses a quantity called the \textit{polynomial capacity} to bound coefficients of certain log-concave polynomials \cite{Gurvits2008}. For a polynomial $p(x_1,\ldots,x_n)$ with nonnegative coefficients and a prescribed exponent vector $\mu$, the corresponding capacity is defined as 
\[
\operatorname{Cap}_{\mu}(p)
=
\inf_{x>0}
\frac{p(x_1,\ldots,x_n)}{x^\mu}.
\]
After the change of variables $x_i=e^{y_i}$ and taking $\log$, the capacity problem is a convex optimization problem.

The goal of the capacity method is to approximate the $x^\mu$ coefficient of a given polynomial $p$ by $\operatorname{Cap}_{\mu}(p)$. This method works particularly well when the polynomial $p$ lies within certain natural log-concave polynomial classes, like real stable, Lorentzian (aka strongly/completely log-concave) \cite{Gurvits2009Newton,BrandenHuh2020Lorentzian,AnariLiuOveisGharanVinzant2024}, and denormalized/dually Lorentzian \cite{brändén2021lowerboundscontingencytables,ross2023duallylorentzianpolynomials}. These classes encode log-concavity inequalities as certain analytic properties of polynomials, and they have had a wide range of recent application in combinatorics and theoretical computer science.

Schur polynomials are denormalized Lorentzian by \cite{huh2019logarithmicconcavityschurrelated} and dually Lorentzian by \cite{ross2023duallylorentzianpolynomials}. Consequently, the capacity method can be applied to obtain bounds on Kostka numbers (their coefficients) using \cite{brändén2021lowerboundscontingencytables}. Although the associated capacity problem is convex and therefore suggests the possibility of efficient computation, a careful algorithmic analysis is still required to establish a polynomial-time approximation algorithm.

The main contribution of this paper is to prove that the capacity can be approximated efficiently for Schur polynomials, thereby providing a framework for obtaining capacity bounds on Kostka numbers. In contrast to related asymptotic and polyhedral results on stretched Kostka coefficients and Kostka polytopes (e.g., see \cite{McAllister2008,Narayanan_2026}), our bounds apply directly to fixed pairs $(\lambda,\mu)$ and are obtained through an efficiently computable convex optimization problem. Our main results are summarized in the following theorem.

\begin{theorem}[Main Theorem]\label{thm:main}
\MainTheoremBody
\end{theorem}

Note that the problem $(P_\mu)$ in the theorem fits into the framework of weighted maximum-entropy dual problems studied in \cite{straszak2019maximumentropydistributionsbit}, and the polynomial dependence on $\log(1/\epsilon)$ is consistent with their results. We note that their theorems are more general and give looser bounds than are possible in our setting; thus we instead use the geometry of the Schur Newton polytope to obtain more specific and sharper results. We also describe how to use their results in \Cref{subsec:bit-complexity-approach}.

This paper is organized as follows. \Cref{sec:prelim} gives background on Schur polynomials, relevant classes of log-concave polynomials, and the capacity method. In \Cref{sec:schur-capacity-bound}, we derive explicit capacity-based bounds for nonzero Kostka numbers. In \Cref{sec:computing-capacity}, we prove that the relevant capacity can be approximated in polynomial time, including a separate treatment of points on faces of the Newton polytope. Finally, in \Cref{sec:summary-examples}, we summarize the results and illustrate the bounds through numerical examples. These examples indicate that the optimization is efficient in practice, and that the quality of the capacity bounds depends on the position of the weight in the Newton polytope.

\subsection{AI Disclosure}
The original mathematical contribution in this paper, including new theorems, lemmas, proofs, and other mathematical arguments, was developed by the authors. AI tools, namely GPT-6 Astra, were used for language polishing, grammar correction, LaTeX typesetting, coding examples, and reviewing mathematical arguments for potential errors. AI tools were also used to draft the abstract. All AI-assisted output was reviewed and verified by the authors.

\subsection{Acknowledgments}

Both authors acknowledge the support of the Natural Sciences and Engineering Research Council of Canada (NSERC), [funding reference number RGPIN-2023-03726]. Cette recherche a \'et\'e partiellement financ\'ee par le Conseil de recherches en sciences naturelles et en g\'enie du Canada (CRSNG), [num\'ero de r\'ef\'erence RGPIN-2023-03726].

\section{Preliminaries}\label{sec:prelim}
\subsection{Notation}\label{subsec:notation}
In the following sections, we let $\mathbb{N} = \{0, 1, ...\}$. For any $\mu=(\mu_1, ..., \mu_n) \in \mathbb{N}^n$, $x^\mu$ denotes $x_1^{\mu_1}...x_n^{\mu_n}$, and $\mu!$ denotes $\prod_{i=1}^n\mu_i!$.  Similarly, for $y=(y_1, ..., y_n) \in \mathbb{R}^n,$ $e^y$ denotes the vector $(e^{y_1}, ..., e^{y_n})$. For $\kappa=(\kappa_1, ..., \kappa_n) \in \mathbb{N}^n$, $\mathbb{R}_\kappa[x_1, ..., x_n]$ denotes the ring of polynomials over $\mathbb{R}$ where the degree of $x_i$ is less than or equal to $\kappa_i$. For a partition $\lambda,$ $l(\lambda)$ denotes its length, $\lambda'$ denotes its conjugate, and $|\lambda|$ denotes its size (the sum of its entries). $||\cdot||_p$ denotes the $p$-vector norm, and $||\cdot||$ denotes the usual 2-norm. $[n]$ denotes the set of integers from 1 to $n$. For a polynomial $p$ and integer vector $\alpha,$ $p_\alpha$ denotes the coefficient of $p$ indexed by the exponent vector $\alpha.$

\subsection{Schur Polynomials}\label{subsec:schur-polynomials}
Schur polynomials are a widely studied class of symmetric functions indexed by partitions. They can be defined as the generating function of semistandard Young tableaux (SSYT) counted by weight. 

\begin{definition}[Schur Polynomial]
    Given a partition $\lambda$, the \textit{Schur polynomial} indexed by $\lambda$ over $n$ variables, where $n \ge l(\lambda)$, is
    \[
    s_\lambda(x_1, ..., x_n)=\sum_{\mu}K_{\lambda, \mu}x^\mu,
    \]
    where $K_{\lambda, \mu}$, called the \textit{Kostka numbers}, is the number of semistandard Young tableaux of size $\lambda$ and weight $\mu$.
\end{definition}

Schur polynomials have been shown to possess a variety of identities. In particular, they are easy to evaluate, as can be seen from the Jacobi--Trudi identity. 

\begin{theorem}[Jacobi-Trudi Identity]\label{thm:jacobi-trudi}
    \[
    s_\lambda(x)=\text{det}[h_{\lambda_i+j-i}(x)]_{i,j=1}^{l(\lambda)}
    \]
    and
    \[
    s_\lambda=\det[e_{\lambda'_i+j-i}]_{i,j=1}^{\lambda_1},
    \]
where $h_k, e_k$ are the complete homogeneous and elementary symmetric functions over $x_1, ..., x_n$, respectively.
\end{theorem}

There is also an explicit formula for $s_\lambda(1^n)$ by the Weyl formula, which is equal to 
$\sum_\mu K_{\lambda, \mu} (1^n)^\mu=\sum_\mu K_{\lambda, \mu}$. 
\begin{theorem}\label{thm:weyl-dimension}
    \[
    s_\lambda(1^n)=\prod_{1 \le i<j\le n}\frac{\lambda_i-\lambda_j+j-i}{j-i}.
    \]
\end{theorem}

Thus, Schur polynomials are easy to evaluate directly. This fact motivates using evaluations of $s_\lambda$, through the capacity method, to bound its coefficients.

\subsection{Log-Concave Polynomials}\label{subsec:log-concave}
Lorentzian, denormalized Lorentzian, and dually Lorentzian polynomials are classes of log-concave polynomials (in the bivariate case, their coefficients form a log-concave sequence after possible transformations). Schur polynomials are denormalized Lorentzian, and many other important polynomial families also fall into Lorentzian-type classes, including several Schur-related families studied in \cite{huh2019logarithmicconcavityschurrelated}.

\begin{definition}[Lorentzian Polynomial, \cite{Gurvits2009Newton,BrandenHuh2020Lorentzian,AnariLiuOveisGharanVinzant2024}]
    Let $p(x)\in \mathbb{R}[x_1, ..., x_n]$ be a degree $d$ homogeneous polynomial. $p$ is called \textit{Lorentzian} if the following conditions hold:\\
    (1) All the coefficients of $p$ are nonnegative\\
    (2) For any $i_1, ..., i_k \in [n]$ and $k \in \mathbb{N}$, $\partial_{i_1}...\partial_{i_k} p$ is either identically zero or log-concave on $\mathbb{R}_{>0}$.
\end{definition}

\begin{definition}[Denormalized Lorentzian Polynomial, \cite{brändén2021lowerboundscontingencytables}] Let $p(x)=\sum_\mu p_\mu x^\mu \in \mathbb{R}[x_1, ..., x_n]$ be a homogeneous polynomial. $p$ is called \textit{denormalized Lorentzian} if its normalization, $N(p(x)):=\sum_\mu \frac{p_\mu}{\mu !} x^\mu$ is Lorentzian.
\end{definition}

\begin{definition}[Dually Lorentzian Polynomial, \cite{ross2023duallylorentzianpolynomials}]
Let $p$ be homogeneous with nonnegative coefficients. We say that $p$ is \textit{dually Lorentzian} if, for some $\kappa$ such that $p\in \mathbb{R}_{\kappa}[x_1,\ldots,x_n]$,
\[
p^\vee(x)=x^\kappa p(x_1^{-1},\ldots,x_n^{-1})
\]
is denormalized Lorentzian.

\end{definition}

\subsection{Capacity Bounds}\label{subsec:capacity-bounds}
The capacity of a polynomial with nonnegative coefficients is defined as an infimum. 
  
\begin{definition}[Capacity]
    For a polynomial $p \in \mathbb{R}_+[x_1, ..., x_n]$ and $\alpha \in \mathbb{N}^n$, the  \textit{capacity} of $p$ is defined as
    \[
    \text{Cap}_\alpha(p(x)):=\text{inf}_{x>0} \frac{p(x)}{x^\alpha}.
    \]

\end{definition}
Capacity can be used to bound the coefficients of log-concave polynomials. An upper bound $p_\alpha \le \text{Cap}_\alpha(p(x))$ is immediate from definition. Theorem 5.10 of \cite{brändén2021lowerboundscontingencytables} gives a lower bound on the coefficients of denormalized Lorentzian polynomials with capacity. Since dually Lorentzian polynomials become denormalized Lorentzian after taking duals, the same method can also be applied to the dual polynomial, giving bounds on the original coefficients after reindexing.

\begin{theorem}[Capacity Bound on Denormalized Lorentzian Polynomials]\label{thm:denorm-capacity}
Let $p \in \mathbb{R}[x_1, ..., x_n]$ be a denormalized Lorentzian polynomial of degree $d$ given by $p(x_1, ..., x_n)=\sum_\mu p_\mu x^\mu$. For any $\alpha \in \mathbb{N}^n$ where $\alpha_1+...+\alpha_n=d$, let $d_i$ be the degree of $x_i$ in $\partial_{i+1}^{\alpha_{i+1}}...\partial_n^{\alpha_n} p|_{x_{i+1}=...=x_n=0}$ for $i=1,2, ..., n-1$ and $d_n$ be the degree of $x_n$ in $p$. Then
\[
\frac{p_\alpha}{\text{Cap}_\alpha(p)}
\ge
\prod_{i=2}^n
\max\left\{
\frac{\alpha_i^{\alpha_i}}{(\alpha_i+1)^{\alpha_i+1}},
\frac{(d_i-\alpha_i)^{d_i-\alpha_i}}{(d_i-\alpha_i+1)^{d_i-\alpha_i+1}}
\right\}.
\]
\end{theorem}

Since $f(x)=\frac{x^x}{(x+1)^{x+1}}$ is a decreasing function, in practice, we can compare the values of $\alpha_i$ and $d_i-\alpha_i$ to give a tighter bound. 
\begin{remark}\label{rem:variable-order}
    Since the proof of \Cref{thm:denorm-capacity} works regardless of the order in which the variables are reduced (see \cite[Theorem 5.10]{brändén2021lowerboundscontingencytables}), for any $\sigma \in S_n$, we can replace the $i$'s in the statement of the theorem by $\sigma(i)$'s. 
\end{remark}

\section{Applying the Capacity Bound on Schur Polynomials}\label{sec:schur-capacity-bound}

The bounds in \Cref{subsec:capacity-bounds} are particularly useful for polynomials satisfying two conditions: (1) they have a suitable log-concavity property, such as being denormalized Lorentzian, and (2) their capacity is easy to evaluate compared to direct coefficient extraction. The discussion in \Cref{sec:prelim} suggests that Schur polynomials may satisfy both conditions: they are denormalized Lorentzian, and since they are easy to evaluate and the capacity optimization problem becomes convex after a change of variables, it is natural to expect that their capacity can also be evaluated efficiently. The following sections show how the capacity method gives bounds on the coefficients of Schur polynomials, a procedure that may also apply to many other polynomials that are proved or conjectured to be denormalized Lorentzian in \cite{huh2019logarithmicconcavityschurrelated}. 

By \cite[Theorem 3]{huh2019logarithmicconcavityschurrelated}, \Cref{thm:denorm-capacity} applies to Schur polynomials. Although the $d_i$s in \Cref{thm:denorm-capacity} can be directly calculated, we relax the bound so that it involves only the capacity, the weight corresponding to the coefficient of interest, and the partition indexing the Schur polynomial. 

For a Schur polynomial $s_\lambda(x_1, ..., x_n)$ indexed by a partition $\lambda=(\lambda_1, ..., \lambda_n)$ with $\lambda_1 \ge...\ge\lambda_n\ge0$, since $s_\lambda$ is symmetric, we only need to consider coefficients of monomials $x^\mu$ where $\mu \in \mathbb{N}^n$ is also a partition and $|\mu|=|\lambda|$. Furthermore, there is a simple condition for $K_{\lambda, \mu}$ to be nonzero (e.g., see \cite[Ex.~2, §2.2]{Fulton1997}) based on dominance orders:

\begin{lemma}\label{lem:kostka-nonzero}
For partitions $\lambda=(\lambda_1,..., \lambda_n), \mu=(\mu_1, ..., \mu_n)$ such that $|\lambda|=|\mu|$, $K_{\lambda, \,\mu}>0$ if and only if, for all $i \in \{1, ..., n\}$,
\[
\sum_{k=1}^i\lambda_k \ge \sum_{k=1}^i \mu_k.
\]
When $\lambda, \mu$  satisfy the relationship above, we say $\lambda$\textit{ dominates} $\mu$, denoted by $\lambda \trianglerighteq\mu$.
\end{lemma}

Thus, we can restrict to the coefficients $x^\mu$ where $\mu$ is a partition, $|\lambda|=|\mu|$ and $\lambda \trianglerighteq\mu$ in the following theorem. 

\begin{theorem}\label{thm:schur-capacity-bound}
    Let $s_\lambda(x_1, ..., x_n)$ be a Schur polynomial indexed by $\lambda$. Let $\mu \in \mathbb{N}^n$ be a partition such that $|\mu|=|\lambda|$ and $\lambda \trianglerighteq\mu$. Then, 
    \[
    \frac{K_{\lambda,\mu}}{\text{Cap}_\mu(s_\lambda(x))}
    \ge
    \prod_{i=1}^{n-1}\max\left\{\frac{\mu_i^{\mu_i}}{(\mu_i+1)^{\mu_i+1}}, \frac{A_i^{A_i}}{(A_i+1)^{A_i+1}}\right\},
    \]
    where $A_i=\sum_{k=1}^i(\lambda_k-\mu_k)$.
\end{theorem}

\begin{proof}
Since $s_\lambda(x)$ is denormalized Lorentzian, by  \Cref{thm:denorm-capacity} and \Cref{rem:variable-order}, \
\[
\frac{K_{\lambda, \mu}}{\text{Cap}_\mu (s_\lambda(x))}
\ge
\prod_{i=1}^{n-1}\max\left\{
\frac{\mu_i^{\mu_i}}{(\mu_i+1)^{\mu_i+1}},
\frac{(d_i-\mu_i)^{d_i-\mu_i}}{(d_i-\mu_i+1)^{(d_i-\mu_i+1)}}
\right\},
\]
where $d_i$ is the degree of $x_i$ in $\partial_1^{\mu_1}...\partial_{i-1}^{\mu_{i-1}}s_\lambda(x)|_{x_1=...=x_{i-1}=0}$ for $2\le i\le n-1$ and $d_1$ is the degree of $x_1$ in $s_\lambda(x)$. 

By definition, $d_i=\max\{\nu_i: K_{\lambda, (\mu_1,..., \mu_{i-1}, \nu_i, \nu_{i+1}, ..., \nu_n)}>0, \text{for some } \nu_{i+1},...\nu_n\}$. By \Cref{lem:kostka-nonzero}, $d_i \le \sum_{k=1}^i\lambda_k-\sum_{k=1}^{i-1}\mu_k$, so $d_i-\mu_i \le A_i$, where $A_i=\sum_{k=1}^i (\lambda_k-\mu_k)$. Since $x^x/(x+1)^{x+1}$ is decreasing on $\mathbb{R}_{\ge0}$, we get \[\frac{K_{\lambda,\mu}}{\text{Cap}_\mu(s_\lambda(x))}\ge\prod_{i=1}^{n-1}\max\{\frac{\mu_i^{\mu_i}}{(\mu_i+1)^{\mu_i+1}}, \frac{A_i^{A_i}}{(A_i+1)^{A_i+1}}\}.\]
\end{proof}

\begin{remark}\label{rem:simple-schur-bound}
Let $\phi(t)=t^t/(t+1)^{t+1}$, with the convention $0^0=1$, and set
\[
m_i=\min\{\mu_i,A_i\},
\] where $A_i$ is as defined in \Cref{thm:schur-capacity-bound}.
Since $\phi(t)$ is decreasing on $\mathbb{R}_{\ge0}$, the product in \Cref{thm:schur-capacity-bound} can be rewritten as
\[
\frac{K_{\lambda,\mu}}{\operatorname{Cap}_\mu(s_\lambda(x))}
\ge
\prod_{i=1}^{n-1}
\frac{m_i^{m_i}}{(m_i+1)^{m_i+1}}.
\]
Moreover,
\[
\phi(t)
=
\frac{1}{t+1}\left(\frac{t}{t+1}\right)^t
\ge
\frac{1}{\mathrm e(t+1)},
\]
because $\left(\frac{t}{t+1}\right)^t$ decreases to $\mathrm e^{-1}$ as $t\to\infty$. Thus, we obtain the simpler but looser bound
\[
K_{\lambda,\mu}
\ge
\operatorname{Cap}_\mu(s_\lambda(x))
\prod_{i=1}^{n-1}\frac{1}{\mathrm e(m_i+1) } \ .
\]
Note that
\[
\prod_{i=1}^{n-1}\frac{1}{e(m_i+1) } \ge \prod_{i=1}^{n-1}\frac{1}{e(\mu_i+1) } \ge \prod_{i=1}^n \frac{1}{e(\mu_i+1)} \ge \left( \frac{n}{e(n+|\lambda|)}\right)^n
\]
by the AM-GM inequality. 
Thus, letting $d = |\lambda|,$ we have 
\[
K_{\lambda,\mu}
\ge
\left(\frac{n}{e(n+d)} \right)^n \operatorname{Cap}_\mu(s_\lambda(x)).\]
\end{remark}

\section{Computing the Capacity of Schur polynomials}\label{sec:computing-capacity}

\Cref{sec:schur-capacity-bound} gives a lower bound on the Kostka numbers of a Schur polynomial in terms of its capacity, the weight of interest, and the partition of the Schur polynomial. This section discusses how this capacity can be computed.

First, recall that $\text{Cap}_\alpha(s_\lambda(x_1, ..., x_n))$ is the optimal value of the optimization problem
\[
\inf_{x>0} \frac{s_\lambda(x)} {x^\alpha}.
\]
After making the substitution $x=e^y$ and taking the logarithm, this can be turned into the unconstrained problem
\[
\inf_y \ln\left(\frac{s_\lambda(e^y)}{e^{\langle y,\alpha\rangle}}\right),
\]
which is equivalent to
\[
\inf_y g_\alpha(y):=\ln \sum_\mu K_{\lambda,\mu}e^{\langle y, \mu-\alpha\rangle}.
\]
The function $g_\alpha(y)$ is convex. Moreover, let $g_\alpha^*$ be the optimal value of this problem; whenever $\alpha \in \text{Newt}(s_\lambda(x))$,
\[
\text{Cap}_\alpha(s_\lambda(x))=e^{g_\alpha^*}
\] (e.g., see \cite[Fact 2.12]{gurvits2020countingmatchingscapacitypreserving}).

Let $d=|\lambda|$. Since $s_\lambda$ is homogeneous of degree $d$, the coefficient indexed by $\alpha$ is zero unless $|\alpha|=d$. By symmetry, when studying nonzero coefficients we may reorder $\alpha$ and assume it is a partition. Under this assumption, we can reduce the feasible region to a linear subspace of $\mathbb{R}^n$: for $y=\gamma1^n+w$, where $\gamma\in \mathbb{R}$ and $w\in \text{Span}\{1^n\}^\perp$, we have
\[
g_\alpha(y)=\ln \sum_\mu K_{\lambda, \mu} e^{\gamma \langle 1^n, \mu-\alpha\rangle+\langle w, \mu-\alpha\rangle}
=\ln \sum_\mu K_{\lambda, \mu} e^{\langle w, \mu-\alpha\rangle}
=g_\alpha(w).
\]

So we can restrict the optimization problem to

\[(P_\alpha) \hspace{1em} \text{inf} \hspace{0.5em} g_\alpha(y)=\ln\sum_\mu K_{\lambda, \mu}e^{\langle y, \mu-\alpha \rangle}\]

\[s.t.\hspace{1em} y_1+...+y_n=0\].

In the following sections, we will continue to refer to an optimizer (if it exists) of $(P_\alpha)$ as $y_\alpha^*$, the optimal value of $(P_\alpha)$ as $g_\alpha^*=\inf g_\alpha$, and $d=|\lambda|$. 

Since $(P_\alpha)$ is convex, it is natural to apply the ellipsoid method. We recall a standard form of the algorithm below (e.g., see \cite{bental_nemirovski_optimization3}).

\begin{theorem}[The Ellipsoid Algorithm]\label{thm:ellipsoid}
   Given any $\beta, R>0$, there is an algorithm which, given a strong first-order counting oracle for a convex function $g$ defined on $\mathbb{R}^n$ and a linear subspace $\mathcal{V}_L \subseteq \mathbb{R}^n$, returns a $y^\circ \in \mathcal{V}_L$ such that 
   \[
   g(y^\circ)\le g^*_R+\beta (g^\bullet_R-g^*_R),
   \]
   where
   \[
   g^*_R=\text{inf}_{||y||_\infty \le R, \hspace{0.2em} y \in \mathcal{V}_L } \hspace{0.3em} g(y),
   \]
   \[
   g^\bullet_R=\text{sup}_{||y||_\infty \le R, \hspace{0.2em} y \in \mathcal{V}_L} \hspace{0.3em} g(y),
   \]
   and the number of calls to the strong first-order counting oracle for $g$ is bounded by a polynomial in $n$, $\log R$, and $\log\frac{1}{\beta}$, where $n$ is the dimension of the ambient vector space. 
\end{theorem}

\subsection{Bounding Box on $||y^*_\alpha||$}\label{subsec:bounding-box}
The ellipsoid algorithm relies on a bounding box on a minimizer of $(P_\alpha)$. This subsection proves such a bound for $\alpha$ in the relative interior of $\text{Newt}(s_\lambda)$. The boundary cases will be handled separately in \Cref{subsec:boundary-case}.

Obtaining this bound requires certain properties of $\text{Newt}(s_\lambda)$. These will be proved using the combinatorial properties of $s_\lambda$ as discussed in \Cref{sec:schur-capacity-bound} as well as the concept of majorization, which can be seen as an extension of the dominance order of partitions in $\mathbb{R}^n$.

\begin{definition}[Majorization]
    For two vectors $u, v \in \mathbb{R}^n$, if \[\sum_{i=1}^nu_i=\sum_{i=1}^n v_i\] and \[\sum_{i=1}^k u_i^\downarrow \ge \sum_{i=1}^k v_i^\downarrow, \] for all $k=1, ..., n-1$, where $u^\downarrow_i, v_i^\downarrow$ denotes the $i^\text{th}$ largest coordinate of $u, v$ respectively, 
    we say that $u$ \textit{majorizes} $v$ and write $u \succ v$.
\end{definition}

We will use the following standard characterization of permutohedra in terms of majorization (e.g., see \cite[4.C.1]{marshall2011majorization}). 

\begin{lemma}[Hardy-Littlewood-P\'olya Theorem on Majorization]\label{lem:hlp-majorization}
    For two vectors $u, v \in \mathbb{R}^n$, \[v \in \text{conv}(\{\sigma(u): \sigma \in S_n\}) \hspace{0.3em} \text{if and only if} \hspace{0.3em} u \succ v\]. 
\end{lemma}

Let $S(\lambda)=\{\sigma(\lambda): \sigma \in S_n\}$. By \Cref{lem:kostka-nonzero,lem:hlp-majorization}, for $z \in \mathbb{N}^n$, we immediately have
\[z \in \text{conv}(S(\lambda)) \iff \lambda  \succ z \iff \lambda \trianglerighteq z^\downarrow \iff z \in \text{supp}(s_\lambda)\]. 

So $\text{supp}(s_\lambda)=\text{conv}(S(\lambda)) \cap\mathbb{N}^n$,  which gives a characterization of $\text{Newt}(s_\lambda)$ as follows. 
\begin{lemma}\label{lem:newton-permutahedron}
\[\text{Newt}(s_\lambda)=\text{conv}(S(\lambda))=\{z\in \mathbb{R}^n: \lambda \succ z\}.\]
\end{lemma}
\begin{proof}
By definition, $\text{Newt} (s_\lambda)=\text{conv}(\text{supp}(s_\lambda))=\text{conv}(\text{conv}(S(\lambda))\cap \mathbb{N}^n).$ So let $C=\text{conv}(S(\lambda)),$ and it suffices to prove that $\text{conv}(C \cap \mathbb{N}^n)=C$. 

Since $\text{Vert}(C)\subseteq S(\lambda) \subset C\cap\mathbb{N}^n, C=\text{conv}(\text{Vert}(C))\subset\text{conv}(C\cap\mathbb{N}^n)$. But $C\cap\mathbb{N}^n \subset C,$ so $\text{conv}(C\cap\mathbb{N}^n)\subset\text{conv}(C)=C$. Thus, $\text{conv}(C\cap\mathbb{N}^n)=C$. 
\end{proof}

Since $\text{Newt}(s_\lambda) \subset\{z \in \mathbb{R}^n: z_1+...+z_n=d\}$, we consider the relative interior instead of the interior of $\text{Newt}(s_\lambda)$. The above lemma immediately gives a simple characterization for being in $\text{relint}(\text{Newt}(s_\lambda))$ whenever the entries of $\lambda$ are not all equal. 

\begin{lemma}\label{lem:relint-sufficient}
Suppose the entries of $\lambda$ are not all equal. Let $z\in \mathbb{R}^n.$ Then, $z \in \text{relint}(\text{Newt}(s_\lambda))$ if and only if 
\[\sum_{i=1}^n z_i=d\] and for $k=1, 2, ..., n-1$, \[\sum_{i=1}^kz_i^\downarrow<\sum_{i=1}^k \lambda_i, \]
. 
\end{lemma}

\begin{lemma}\label{lem:support-halfspace}
    If $z\in \text{Newt}(s_\lambda)$, then $\forall y\in \mathbb{R}^n$, $\text{supp}(s_\lambda) \cap\{w \in \mathbb{R}^n: \langle w-z, y\rangle\ge0\} \neq \emptyset$. 
\end{lemma}

\begin{proof}
By definition, $z = \sum_{j=1}^m \gamma_j w_j,$ where $\sum_{j=1}^m \gamma_j =1, \gamma_1, ..., \gamma_m >0$, $w_1, ..., w_m \in \text{supp}(s_\lambda)$. Then, $0 = \langle z-z, y\rangle = \langle \sum_{j=1}^m \gamma_j (w_j-z), y \rangle = \sum_{j=1}^m \gamma_j \langle  (w_j-z), y\rangle.$ Thus, $\exists l \in \{1, ..., m\}$ such that $\langle w_l-z, y\rangle \ge 0 $
\end{proof}

With this lemma, we can prove a bounding box for $y^*_\alpha$ where the partition $\alpha$ satisfies a more restricted version of the conditions in \Cref{lem:kostka-nonzero}, or, in other words, when $\alpha$ is in the relative interior of $\text{Newt}(s_\lambda)=\text{conv}(S(\lambda))$. 

First, we need to prove that points around $\alpha$ are also in $ \text{Newt}(s_\lambda)$. 

\begin{lemma}\label{lem:relint-perturbation}
For a partition $\alpha \in \mathbb{N}^n$, if \[\sum_{i=1}^n \alpha_i=d\] and for all $k=1, 2, ..., n-1$, \[\sum_{i=1}^k \alpha_i <\sum_{i=1}^k \lambda_i,\] then $\forall y\in \mathbb{R}^n$ such that $\sum_{i=1}^n y_i=0$ and $y\ne0$, \[\alpha+\frac{1}{n}\frac{y}{||y||} \in \text{Newt}(s_\lambda).\]
\end{lemma}
\begin{proof}
    Let $z=\alpha+\frac{1}{n}\frac{y}{||y||}.$ Clearly, $\sum_{i=1}^nz_i=\sum_{i=1}^n\alpha_i=d$. Write $z^\downarrow=(z^\downarrow_1, ..., z_n^\downarrow)$, where $z_1^\downarrow\ge...\ge z_n^\downarrow$. We then have $z_i^\downarrow \le \alpha_i+\frac{1}{n}$. Thus, for $k=1, 2, ..., n-1, \sum_{i=1}^kz^\downarrow_i \le \sum_{i=1}^k(\alpha_i+\frac{1}{n})<1+\sum_{i=1}^k\alpha_i\le \sum_{i=1}^k\lambda_i$ since both $\sum_{i=1}^k\alpha_i, \sum_{i=1}^k \lambda_i$ are integers. By \Cref{lem:newton-permutahedron}, $z \in \text{Newt}(s_\lambda)$. 
\end{proof}

Using the above lemma, we can finally prove the bounding box on $y^*_\alpha$.
\begin{theorem}\label{thm:optimizer-box}
For a partition $\alpha \in \mathbb{N}^n$ such that $\alpha \in \text{relint}(\text{Newt}(s_\lambda))$, there exists an optimizer $y^*_\alpha$ of $(P_\alpha)$ that satisfies \[||y^*_\alpha|| \le n\ln\sum_\mu K_{\lambda, \mu}. \]
\end{theorem}
\begin{proof}

We may assume that $\lambda$ does not have all-equal entries because otherwise $g_\alpha(y) \equiv 0$ and the theorem automatically holds. Thus, $\alpha$ satisfies the hypothesis in \Cref{lem:relint-perturbation}.

Note that $K_{\lambda, \mu}\ge 1$ for any $\mu$ where $K_{\lambda, \mu} \neq 0$ since $K_{\lambda, \mu}$ is an integer.
In the proof, we extend $K_{\lambda, \mu}$ to a function on $\mathbb{R}^n$ (i.e., define $K_{\lambda, \mu}=0$ when $\mu$ is not in $\mathbb{N}^n$). 

Let $0 \neq y \in \mathbb{R}^n$ such that $y_1+...+y_n =0$. Let $z=\alpha+\frac{1}{n}\frac{y}{||y||}.$ By \Cref{lem:relint-perturbation}, $z \in \text{Newt}(s_\lambda)$.  Write $g_\alpha(y)$ as \[\ln(\sum_{ \langle \mu-z, y\rangle \ge0}K_{\lambda, \mu} e^{\langle y, \mu-\alpha \rangle}+\sum_{ \langle \mu-z, y\rangle <0}K_{\lambda, \mu} e^{\langle y, \mu-\alpha \rangle}) \ge \ln\sum_{\langle \mu-z, y\rangle \ge0}K_{\lambda, \mu} e^{\langle y, \mu-\alpha \rangle}.\]

Since $\langle\mu-z, y\rangle=\langle \mu-\alpha-\frac{1}{n}\frac{y}{||y||}, y\rangle$, $\langle\mu-z,y\rangle \ge0 \iff \langle\mu-\alpha, y\rangle \ge \frac{1}{n}||y||$. Since $z \in \text{Newt}(s_\lambda)$,   by \Cref{lem:support-halfspace}, $\exists \mu \in \text{supp}(s_\lambda)$ such that $\langle \mu-z, y\rangle\ge 0$. Thus, 

\[g_\alpha(y)\ge \ln\sum_{ \langle \mu-z, y\rangle \ge0}K_{\lambda, \mu} e^{\langle y, \mu-\alpha \rangle} \ge \ln(e^{\frac{1}{n}||y||})=\frac{1}{n}||y||.\]

This shows that $g_\alpha$ is coercive so there exists an optimizer $y^*_\alpha$ on the linear subspace $\{y \in \mathbb{R}^n: y_1+...+y_n =0\}$. We may assume that $y_\alpha^* \neq 0$ because otherwise the theorem also automatically holds. 

Plugging in $y=0^n$ gives an upper bound on $g(y^*_\alpha)$, 
\[g_\alpha(y^*_\alpha)\le g_\alpha(0^n)=\ln\sum_\mu K_{\lambda, \mu}. \]

Rearranging, we obtain 
\[||y^*_\alpha||\le n\ln\sum_\mu K_{\lambda, \mu}.\]
\end{proof}

This bound on a minimizer of $(P_\alpha)$ involves the parameter ln($\sum_{\mu} K_{\lambda, \mu}$), which can be further bounded as follows using \Cref{thm:weyl-dimension}.  

\begin{lemma}\label{lem:weyl-log-bound}
    \[\ln\sum_\mu K_{\lambda, \mu} \le (\ln(n)+\gamma+\frac{1}{2n})\sum_{i=1}^{n-1}\lambda_i=(\ln(n)+\gamma+\frac{1}{2n})(d-\lambda_n), \]
    where $\gamma\approx0.5772$ is the Euler–Mascheroni constant.  
\end{lemma}
\begin{proof}
By \Cref{thm:weyl-dimension}, ln$\sum_{\mu} K_{\lambda, \mu}=\ln (s_\lambda(1^n))=\sum_{1\le i<j \le n} \ln(1+\frac{\lambda_i-\lambda_j}{j-i})$. Since $\ln(1+0)=0, \frac{d}{dx}\ln(1+x)=\frac{1}{1+x}\le1=\frac{d}{dx} x, \forall x\ge 0, $ we have $\ln(1+x)\le x, \forall x\ge 0$. Thus, 
\[\ln(\sum_\mu K_{\lambda, \mu}) \le \sum_{1\le i<j\le n}\frac{\lambda_i-\lambda_j}{j-i} \le\sum_{1\le i<j\le n}\frac{\lambda_i}{j-i}=\sum_{i=1}^{n-1} \lambda_i\sum_{k=1}^{n-i} \frac{1}{k} \].

  $\sum_{k=1}^m\frac{1}{k}\le \ln(m)+\gamma+\frac{1}{2m}$, so by replacing the $\sum_{k=1}^{n-i} \frac{1}{k}$ of each $i$ by the approximation for the harmonic sums from 1 to $n$,  we get the result. 
\end{proof}

Let $\alpha \in \text{relint}(\text{Newt}(s_\lambda)), \mathcal{V}_L=\{y \in \mathbb{R}^n: y_1+...+y_n=0\}, R=\max\{1, n\ln\sum_\mu K_{\lambda, \mu}\}$. By \Cref{thm:optimizer-box} and \Cref{lem:weyl-log-bound}, the set $\{y\in \mathbb{R}^n: ||y||_\infty \le R\}$ contains a minimizer of $(P_\alpha)$. Thus, $(P_\alpha)$ can be rewritten as the problem ($P_\alpha'$), which is the one that we consider in \Cref{subsec:counting-oracles,subsec:relint-algorithm}:

\[(P_\alpha') \hspace{1em}\text{inf}\hspace{0.3em} g_\alpha(y)=\ln\sum_\mu K_{\lambda, \mu} e^{\langle y, \mu-\alpha \rangle}\]
\[s.t. \hspace{0.3em} y \in \mathcal{V}_L,\]
\[||y||_\infty \le R.\]

\subsection{Counting Oracles for the Ellipsoid Algorithm}\label{subsec:counting-oracles}
In addition to a bounding box on the optimizer, the ellipsoid algorithm requires first-order information for $g_\alpha(y)$. In this subsection, we show that both the objective function and its gradient can be approximated efficiently.

More precisely, we show that, given an $\epsilon >0$ and $y \in \mathbb{R}^n$ such that  $||y||_\infty \le R$, there exists an oracle that can return a $\tilde{g}_\alpha(y)\in \mathbb{R}$ such that $|g_\alpha(y)-\tilde{g}_\alpha(y)|\le \epsilon$ and a $\tilde{G}(y) \in \mathbb{R}^n$ such that $||\tilde{G}(y)-\nabla g_\alpha(y)||_1 \le \epsilon$ that is polynomial-time in a set of parameters.  To see this, we rewrite the objective function as
$g_\alpha(y)=\ln (s_\lambda(e^{y}))-\langle y, \alpha\rangle$, use the Jacobi-Trudi identity to evaluate $s_\lambda$ and its logarithmic derivatives, and then prove Lipschitz bounds that control the error caused by approximating the exponential and logarithm.

First, \Cref{lem:h-computation,lem:schur-computation,lem:schur-log-gradient-computation} give the bit complexities of evaluating the Schur polynomial and its log-gradient on rational inputs with prescribed bit-lengths using the Jacoobi-Trudi identity.

\begin{lemma}\label{lem:h-computation}
Let $h_k(x_1, ..., x_n)$ be the complete homogeneous symmetric polynomial on $n$ variables of degree $k$, where $k \in \mathbb{N}$. Suppose $x\in(\mathbb{Q}_{>0})^n$ and the bit-length of each entry of $x$ is bounded by $M$, then $h_k(x)$ can be computed in time polynomial in $n$, $k$, and $M$.
\end{lemma}

\begin{proof}
Letting $h_i(x_1)=x_1^i, \forall i \in \mathbb{N} $ and $h_0(x_1)=...=h_0(x_1, ..., x_n)=1$, we have the recurrence relationship that $h_i(x_1, ..., x_j)=x_jh_{i-1}(x_1, ..., x_j)+h_i(x_1, ..., x_{j-1}), \forall i\ge1, j\ge2$. Thus, we can compute $h_k(x_1, ..., x_n)$ with the table
\[
\begin{array}{c|cccc}
 & 0&1 &   \cdots & k \\
\hline
(x_1) & 1 & x_1 & \cdots & x_1^k \\
(x_1,x_2) & 1 & h_1(x_1,x_2) & \cdots & h_k(x_1,x_2) \\
\vdots & \vdots & \vdots & \ddots & \vdots \\
(x_1,\dots,x_n) & 1& h_1(x_1,\dots,x_n) & \cdots & h_k(x_1,\dots,x_n)
\end{array}
\]
where each $h_i(x_1, ..., x_j)$ can be calculated as the sum of $x_j$ multiplied by the entry $h_{i-1}(x_1, ..., x_j)$ on its left with the entry $h_i(x_1, ..., x_{j-1})$ above it, so the total number of calculations is $nk$. 

Since every table entry $h_i(x_1,\dots,x_j)$ is a sum of
$\binom{i+j-1}{i}$ monomials, each a product of $i$ rationals whose
bit-lengths are at most $M$, its bit-length is bounded by
\[
O(ijM+\log \binom{i+j-1}{i})
\le O(nkM+k\log(n+k)).
\]
Thus all intermediate values have polynomial bit-length in $n,k,M$.
There are $O(nk)$ table entries, and each entry is computed using one
multiplication and one addition of rationals of polynomial bit-length.
Therefore the total running time is polynomial in $n,k,M$.
\end{proof}

\begin{lemma}\label{lem:schur-computation}
Let $s_\lambda(x_1, ..., x_n)$ be a Schur polynomial on $n$ variables indexed by the partition $\lambda=(\lambda_1, ..., \lambda_n)$ of length $l$. Suppose $x\in (\mathbb{Q}_{>0})^n$ and the bit-length of each entry of $x$ is bounded by $M$, then $s_\lambda(x)$ can be computed in time polynomial in $n$, $\lambda_1$, and $M$, and the output also has bit-length polynomial in these parameters. 
\end{lemma}

\begin{proof}
By the Jacobi--Trudi identity,
$s_\lambda(x)=\det(H)$, where  
$H=[h_{\lambda_i+j-i}(x)]_{i,j=1}^l$ and we use the convention $h_0=1$ and $h_r=0$ for $r<0$.

Since $\lambda_i+j-i \le \lambda_1+l-1 \le \lambda_1+n-1$, 
\Cref{lem:h-computation} implies that all entries of $H$ can be computed in time
polynomial in $n,\lambda_1,M$, and their bit-lengths are also polynomial
in these parameters.

Finally, the determinant of the $l\times l$ matrix $H$ can be computed
by Gaussian elimination using arithmetic on rationals of polynomial
bit-length. Since $l\le n$, this takes polynomial time in
$n,\lambda_1,M$. Therefore $s_\lambda(x)$ can be computed in
$\operatorname{poly}(n,\lambda_1,M)$ time, and the output also has bit-length polynomial in these parameters.
\end{proof}

\begin{lemma}\label{lem:schur-log-gradient-computation}
  Let $s_\lambda(x_1, ..., x_n)$ be a Schur polynomial on $n$ variables indexed by the partition $\lambda=(\lambda_1, ..., \lambda_n)$ of length $l$. If $x\in (\mathbb{Q}_{>0})^n$ and the bit-length of each entry of $x$ is bounded by $M$, the partial derivatives $\partial_i\ln(s_\lambda(x)), \forall i =1, ..., n$, can be computed in time polynomial in $n, \lambda_1,$ and $ M$. 
\end{lemma}
\begin{proof}
Since $\partial_i h_j(x_1,\dots,x_n)
=
\sum_{r=0}^{j-1}(r+1)x_i^r
h_{j-r-1}(x_1,\dots,x_{i-1},x_{i+1},\dots,x_n)$
and the degrees satisfy
$j\le \lambda_1+l-1\le \lambda_1+n-1$ for all complete homogeneous
symmetric polynomials appearing in $H$, the same argument as in the proof
of \Cref{lem:schur-computation} shows that the entries of $\partial_iH$ can be computed in
time polynomial in $n$, $\lambda_1$, and $M$, and have bit-length
polynomial in these parameters.

Since $x\in(\mathbb Q_{>0})^n$, we have $s_\lambda(x)>0$, so $H$ is
invertible. Moreover, $H$ is an $l\times l$ rational matrix with
$l\le n$ whose entries have polynomial bit-length. Therefore, by the
polynomial-time complexity of Gaussian elimination over $\mathbb Q$,
$H^{-1}$ can be computed in time polynomial in $n$, $\lambda_1$, and
$M$. Consequently, the matrix product
$H^{-1}\partial_iH$ and its trace can also be computed in polynomial
time. Thus, $\partial_i\log(s_\lambda(x))$ can be computed in time
polynomial in $n$, $\lambda_1$, and $M$.
\end{proof}

Next, we prove Lipschitz bounds for the function $f(x)=\ln(s_\lambda(e^x))$ and its gradient. 

\begin{lemma}\label{lem:log-schur-lipschitz}
    Let $x, y \in \mathbb{R}^n$. Let  $\delta =||x-y||_\infty$. Then
    \[
    |\ln(s_\lambda(e^x))-\ln(s_\lambda(e^y))| \le \delta d,
    \]
    where $d=|\lambda|$. 
\end{lemma}
\begin{proof}
 The inequality can be obtained from the following calculations:
\[\ln(s_\lambda(e^x))-\ln(s_\lambda(e^y))=\ln\frac{s_\lambda(e^x)}{s_\lambda(e^y)}=\ln \sum_\mu \frac{K_{\lambda, \mu} e^{\langle x, \mu \rangle}}{\sum_\mu K_{\lambda, \mu}e^{\langle y, \mu \rangle}}\] \[=\ln \sum_\mu \left(\frac{K_{\lambda, \mu} e^{\langle y, \mu\rangle}}{\sum_\mu K_{\lambda, \mu}e^{\langle y, \mu \rangle}}\cdot\frac{e^{\langle x, \mu \rangle}}{e^{\langle y, \mu \rangle}}\right)\le \ln \max_\mu \frac{e^{\langle x, \mu \rangle}}{e^{\langle y, \mu \rangle}} \]\[\le \max_\mu \langle x-y, \mu\rangle \le ||x-y||_\infty||\mu||_1=\delta d.\] 
Similarly, $\ln(s_\lambda(e^y))-\ln(s_\lambda(e^x)) \le \delta d$, so the result follows.
\end{proof}

\begin{lemma}\label{lem:gradient-lipschitz}
    Let $x, y \in \mathbb{R}^n$. Let $ \delta = ||x-y||_\infty$. Then, 
    \[||\nabla\ln(s_\lambda(e^x))-\nabla\ln(s_\lambda(e^y))||_1 \le 2\delta d^2 , \]
    where $d=|\lambda|$. 
\end{lemma}

\begin{proof}
Let $f(x)=\ln(s_\lambda(e^x))$. First,
\[
\nabla f(x)
=\frac{\sum_\mu K_{\lambda, \mu}e^{\langle x, \mu \rangle}\mu }{s_\lambda(e^x)}
=\frac{\sum_\mu K_{\lambda, \mu} e^{\langle x, \mu\rangle}\mu}{\sum_\mu K_{\lambda, \mu} e^{\langle x, \mu \rangle}}
=\sum_\mu w_\mu (x) \mu,
\]
where
\[
w_\mu (x):=\frac{K_{\lambda, \mu}e^{\langle x, \mu \rangle}}{\sum_\mu K_{\lambda, \mu} e^{\langle x,\mu \rangle}} \ge 0,
\qquad
\sum_\mu w_\mu(x)=1.
\]
Also, with this notation, 
\[\nabla^2 f(x)=\sum_\mu w_\mu(x) \mu\mu^T-(\sum_\mu w_\mu(x) \mu)(\sum_\mu w_\mu(x)\mu)^T.\]
Denoting the $j^{\text{th}}$ row of $\nabla^2 f(x)$ by $(\nabla^2 f(x))_{j}$ and the $(i,j)$-entry of $\nabla^2 f(x)$ by $(\nabla^2 f(x))_{i,j}$, we have that 
\[(\nabla^2 f(x))_{i,j}=\sum_\mu w_\mu(x)\mu_i\mu_j-(\sum_\mu w_\mu(x)\mu)_i(\sum_\mu w_\mu(x)\mu)_j, \]
and
\[\sum_{i,j}|(\nabla^2 f(x))_{i,j}| \le \sum_{i,j} \sum_\mu w_\mu(x) \mu_i\mu_j+\sum_{i,j}(\sum_\mu w_\mu(x)\mu)_i(\sum_\mu w_\mu(x)\mu)_j\]\[=\sum_\mu w_\mu(x)\sum_{i,j} \mu_i\mu_j+\sum_i(\sum_\mu w_\mu(x) \mu)_i\sum_j (\sum_\mu w_\mu(x)\mu)_j \le \sum_\mu w_\mu(x) d^2+d^2=2d^2. \]
By the FTC, $\nabla f(x)-\nabla f(y)= \int_0^1 \nabla^2f(y+t(x-y))(x-y)dt$, so taking $||\cdot||_1$s on both sides, 
\[||\nabla f(x)-\nabla f(y)||_1 =||\int_0^1\nabla^2 f(y+t(x-y))(x-y)dt||_1\]
\[\le \int_0^1||\nabla^2 f(y+t(x-y))(x-y)||_1dt  =\int_0^1 \sum_{j=1}^n |\langle (\nabla^2 f(y+t(x-y))_j, x-y\rangle|dt\]
\[\le ||x-y||_\infty \int_0^1\sum_{j=1} ^n||(\nabla ^2 f(y+t(x-y)))_j||_1dt \]
\[=||x-y||_\infty \int_0^1 \sum_{i,j} |(\nabla^2 f(y+t(x-y)))_{i,j}|\le 2\delta d^2.\]
\end{proof}

Since we need to approximate $e^y$ with a rational, we will make use of the following lemma. 

\begin{lemma}\label{lem:log-approx}
For $w\in\mathbb{R}_{>0}, u \in \mathbb{R}$,  if $\frac{|e^u-w|}{|e^u|} < 2^{-B-1}$ for some $B\in \mathbb{N}$, then $|\log (w)-u|<2^{-B}.$
\end{lemma}

\begin{proof}
    \[|\log(w)-u|=|\log(\frac{w}{e^u})|=|\log(1+\frac{w-e^u}{e^u})|.\]
    Since $|\frac{e^u-w}{e^u}|<1/2$, \[1+\frac{w-e^u}{e^u} \in (\frac{1}{2}, \frac{3}{2}).\] Since for all $t$ on $(1/2, 3/2), \frac{d}{dt}\log(t)\le 2,$
    \[|\log(1+\frac{w-e^u}{e^u})|= |\log(1+\frac{w-e^u}{e^u})-\log(1)| \le 2|\frac{w-e^u}{e^u}|<2^{-B}\] by the Mean Value Theorem.
\end{proof}

 With these lemmas, we are ready to prove the existence of an approximate counting oracle for $g_\alpha$.

\begin{theorem}[Existence of an Approximate Counting Oracle for $g_\alpha(y)$]\label{thm:g-oracle}
There exists an approximate counting oracle that, given $0<\epsilon<1$ and $y \in \mathbb{R}^n$ such that $||y||_\infty \le R$, returns a $\tilde{g}_\alpha(y)$ such that \[g_\alpha(y)-\epsilon \le \tilde{g}_\alpha(y) \le g_\alpha(y)+\epsilon,\] with running time polynomial in $n, d$, and $\log(\frac{1}{\epsilon})$.
\end{theorem}
\begin{proof}
By \Cref{lem:log-schur-lipschitz}, if we choose a rational approximation $y'$ of $y$ with bit precision $B>\log_2(\frac{4d}{\epsilon})$, $||y'-y||_\infty \le \frac{\epsilon}{4d},$ so $|\ln(s_\lambda(e^y))-\ln(s_\lambda(e^{y'}))|<\frac{\epsilon}{4}$. Since $||y||_\infty \le R,$ $||y'||_\infty\le R+\frac{\epsilon}{4d}<R+\frac{1}{4}$. Also, the error produced by using $y'$ to approximate $\langle y, \alpha \rangle $ is bounded by $\frac{\epsilon}{4}$, since $|\alpha| = d.$

Next, we want to obtain an approximation $z_i$ to each $e^{y'_i}$ such that $|\log(z_i)-y'_i| < \frac{\epsilon}{4d}$ so that $|\ln(s_\lambda(z))-\ln(s_\lambda(e^{y'}))|<\frac{\epsilon}{4}.$ By \Cref{lem:log-approx}, the bound on $||y'||_\infty, $ and standard polynomial-time approximation algorithms for the exponential function (e.g., see \cite{Brent1976}), we can find  a rational approximation $z_i$ of $e^{y'_i}$ with relative error  within $\frac{\epsilon}{8d}$ in time polynomial in $\log(\frac{1}{\epsilon})$, $\log(d),$ and $R$. Each $z_i$ would then have bit-length polynomial in the same parameters. 

Then, by \Cref{lem:schur-computation}, $s_\lambda(z)$ can be computed in time polynomial in these parameters as well as $\lambda_1$ and $n$, with the output also having bit-length polynomial in them.

Since $||\log(z)-y'||_\infty < \frac{\epsilon}{4d}<\frac{1}{4d}$, $||y-y'||_\infty <\frac{1}{4d}$, and $||y||_\infty \le R$,
\[
e^{\langle\mu, \log(z)\rangle}
\ge e^{-||\mu||_1||\log(z)||_\infty}
> e^{-dR-1/2}=e^{-dR}/e^{1/2},
\]
for every exponent vector $\mu$ appearing in $s_\lambda$. Thus, $s_\lambda(z)\ge e^{-dR}/e^{1/2}>0$, as nonzero coefficients of the Schur polynomial are positive integers and there is at least one monomial. Thus, using standard polynomial-time approximation algorithms for the logarithm (e.g., see \cite{Brent1976}), we can approximate $\log(s_\lambda(z))$ with error bounded by $\frac{\epsilon}{4}$ with time polynomial in $d, R, \log(\frac{1}{\epsilon})$, $\lambda_1$ and $n$. The total approximation error is bounded by $\frac{\epsilon}{4}\cdot 4=\epsilon$. By \Cref{lem:weyl-log-bound}, $R$ is polynomial in $n$ and $d$. Using the fact that $\lambda_1 \le d$, we have that the time complexity is polynomial in $n, d,$ and $\log(\frac{1}{\epsilon})$.
\end{proof}

Using \Cref{lem:schur-log-gradient-computation}, \Cref{lem:gradient-lipschitz}, and \Cref{lem:log-approx} the existence of an approximate first-order oracle for $\nabla g_\alpha(y)$ follows by an argument analogous to the proof of \Cref{thm:g-oracle}.

\begin{theorem}[Existence of an Approximate Counting Oracle for $\nabla g_\alpha(y)$]\label{thm:grad-oracle}
There exists an approximate counting oracle that, given $0<\epsilon<1$ and $y \in \mathbb{R}^n$ such that $||y||_\infty \le R$, returns a $\tilde{G}(y)$ such that \[||\nabla g_\alpha(y)-\tilde{G}(y)||_1 \le \epsilon,\] with running time polynomial in $n, d$, and $\log(\frac{1}{\epsilon})$. 
    
\end{theorem}

\subsection{Algorithm Time Complexity for \(\alpha \in \text{relint}(\text{Newt}(s_\lambda))\) }\label{subsec:relint-algorithm}
Given the bounding box from \Cref{subsec:bounding-box} and the approximate first-order oracle from \Cref{subsec:counting-oracles}, the ellipsoid algorithm can be applied to approximate $\text{Cap}_\alpha(s_\lambda)$. This subsection analyzes the number of oracle calls and then derives the time complexity for $(P_\alpha)$ when $\alpha \in \text{relint}(\text{Newt}(s_\lambda))$. The standard robustness theory for the ellipsoid method with approximate first-order oracles applies in our setting, following the framework of \cite[Section 6.2]{singh2013entropyoptimizationcounting}. Specifically, their argument applies to $g_\alpha(y)$ after replacing their $\sqrt{m}$ by $2d=2|\lambda|$ and replacing their stopping radius by $\epsilon/(32d)$. Thus, we address two scenarios: (1) there exists an exact first-order counting oracle for $g_\alpha(y)$ and (2) there only exists an approximate first-order counting oracle for $g_\alpha(y)$.

\begin{theorem}\label{thm:ellipsoid-exact-oracle}
Given a partition $\alpha \in \text{relint}(\text{Newt}(s_\lambda))$, for every $\epsilon>0$, there exists an algorithm that, given a strong first-order counting oracle for $g_\alpha(y)$, returns a $y^\circ \in \mathbb{R}^n$ such that
\[
g_\alpha(y^\circ) \le g^*_\alpha+\epsilon,
\]
where $g_\alpha^*$ is the optimal value of $(P_\alpha)$, and the number of calls to the strong first-order counting oracle is polynomial in $n, \log(d)$, and $\log(\frac{1}{\epsilon})$.
    
\end{theorem}
\begin{proof}
For any feasible $y$ to $(P_\alpha')$, we have \[g_\alpha(y) \le \ln\sum_\mu K_{\lambda, \mu} e^{||y||_\infty ||\mu-\alpha||_1}\le \ln\sum_{\mu}  K_{\lambda, \mu} e^{2dR}=2dR+\ln\sum_\mu K_{\lambda, \mu}\] and 
 \[g_\alpha(y) \ge \ln\sum_\mu K_{\lambda, \mu} e^{-||y||_\infty||\mu-\alpha||_1} \ge \ln\sum_\mu K_{\lambda, \mu} e^{-2dR}=-2dR+\ln \sum_\mu K_{\lambda, \mu},\]
  using the fact that $\alpha, \mu$ are weak compositions both of size $d$.

Thus, \[(g_\alpha)_R^\bullet:=\text{sup}_{||y||_\infty \le R, \hspace{0.2em}y \in \mathcal{V}_L} g_\alpha(y)\le 2dR+\ln\sum_\mu K_{\lambda, \mu}, \] 

\[(g_\alpha)^*_R:=\text{inf}_{||y||_\infty \le R, \hspace{0.2em} y \in \mathcal{V}_L} g_\alpha(y) \ge -2dR+\ln\sum_\mu K_{\lambda, \mu},\] 
and \[(g_\alpha)^\bullet_R-(g_\alpha)^*_R\le 4dR.\]

Let $\beta=\frac{\epsilon}{4dR}$. By \Cref{thm:ellipsoid}, there exists an algorithm that returns a $y^\circ \in \mathcal{V}_L$ such that \[g_\alpha(y^\circ) \le (g_\alpha)^*_R+\beta((g_\alpha)^\bullet_R-(g_\alpha)^*_R) \le (g_\alpha)^*_R+\beta \cdot 4dR=(g_\alpha)^*_R+\epsilon=g^*_\alpha+\epsilon,\] since, as argued in \Cref{subsec:bounding-box}, ($P_\alpha$) and ($P_\alpha'$) are equivalent as optimization problems, and the number of calls to the strong first-order counting oracle is bounded by a polynomial in $n, \log(R)$, and $ \log(\frac{1}{\beta})=\log(\frac{4dR}{\epsilon})$.

Since $R=O(dn\log(n))$ by \Cref{lem:weyl-log-bound}, we have that $\log(\frac{4dR}{\epsilon})=O(\log(\frac{1}{\epsilon})+\log(n)+\log(d))$. Thus, the number of calls to the strong first-order counting oracle is polynomial in $n$, $\log(d)$, and $\log(\frac{1}{\epsilon})$.
\end{proof}

From the proof of the previous theorem and the underlying ellipsoid algorithm with an approximate first-order counting oracle described in \cite[Section 6.2]{singh2013entropyoptimizationcounting}, we obtain the following theorem.

\begin{theorem}\label{thm:ellipsoid-approx-oracle} 
Given a partition $\alpha \in \text{relint}(\text{Newt}(s_\lambda))$, for every $\epsilon>0$, there exists an algorithm that, given an approximate first-order counting oracle for $g_\alpha(y),$  returns a $y^\circ \in \mathbb{R}^n$ such that
\[
g_\alpha(y^\circ) \le g^*_\alpha+\epsilon,
\]
where $g_\alpha^*$ is the optimal value of $(P_\alpha)$, and the number of calls to the approximate first-order counting oracle is polynomial in $n, \log(d),$ and $\log(\frac{1}{\epsilon})$.
    
\end{theorem}

Finally, using the approximate first-order counting oracles from \Cref{subsec:counting-oracles} in the ellipsoid algorithm, we obtain a polynomial-time algorithm to compute the capacity of Schur polynomials, which is summarized in the following theorem. 

\begin{theorem}\label{thm:relint-capacity-algorithm}
Let $s_\lambda(x) \in \mathbb{R}[x_1, ..., x_n]$ be a Schur polynomial on $n$ variables, and let $d=|\lambda|$. 
For any partition $\alpha \in \text{relint}(\text{Newt}(s_\lambda))$, there exists an optimization problem $(P_\alpha)$ with objective function $g_\alpha(y)$ and an algorithm that, given $\epsilon>0$, returns a $y^\circ \in \mathbb{R}^n$ such that
\[
g_\alpha(y^\circ) \le g^*_\alpha+\epsilon,
\]
where $g^*_\alpha$ is the optimal value of $(P_\alpha)$, and
\[
\text{Cap}_\alpha(s_\lambda(x))=e^{g^*_\alpha}.
\]
The running time of the algorithm is polynomial in $n, d$, and $\log(\frac{1}{\epsilon}).$
    
\end{theorem}

Here, $g_\alpha(y)$ and $(P_\alpha)$ are as defined in \Cref{sec:computing-capacity}. The running time is obtained by directly applying \Cref{thm:g-oracle,thm:grad-oracle} (complexities of our approximate counting oracles for $g_\alpha(y)$ and $\nabla g_\alpha(y)$, respectively), and \Cref{thm:ellipsoid-approx-oracle} (number of calls to an approximate first-order counting oracle for $g_\alpha(y)$ in the underlying ellipsoid algorithm). 

\subsection{Boundary Case}\label{subsec:boundary-case}

This subsection addresses the case where $\alpha$ lies on the boundary of $\text{Newt}(s_\lambda)$. We may assume that the entries of $\lambda$ are not all equal, since that case is already analyzed in the previous sections. As usual, we assume that $\alpha$ is a partition dominated by $\lambda$. The boundary case then occurs exactly when at least one proper dominance inequality is tight; that is, there exists some $j<n$ such that $\alpha_1+...+\alpha_j = \lambda_1+...+\lambda_j.$

Assume that $\alpha$ satisfies the properties described above. First, we show that $\text{Cap}_\alpha(s_\lambda) = \text{Cap}_\alpha((s_\lambda)_F)$, where 
\[
s_\lambda(x)_F := \sum_{\mu \in F} K_{\lambda, \mu} x^\mu
\]
is the Schur polynomial restricted to the face $F$, and $\alpha \in F$. 

\begin{lemma}\label{lem:face-capacity}
    Let $p(x)=\sum_\beta c_\beta x^\beta \in \mathbb{R}_+[x_1, ...,x_n]$. Let $F$ be a face of $\text{Newt}(p)$ where $\alpha \in F$. Then, $\text{Cap}_\alpha(p)= \text{Cap}_\alpha(p_F),$ where $p_F(x) = \sum_{\beta\in F} c_\beta x^\beta.$
\end{lemma}

\begin{proof}
It is immediate that $\text{Cap}_{\alpha} (p) \ge \text{Cap}_\alpha (p_F)$ as $p$ has nonnegative coefficients. 

To show that $\text{Cap}_{\alpha}(p) \le \text{Cap}_\alpha (p_F),$ let $u$ be a normal vector to $F$, i.e.,  $\langle u, \beta\rangle =m, \forall \beta \in F$ and $\langle u, \beta \rangle < m, \forall \beta \in \text{Newt}(p)\setminus F$. For every $x>0,$
\[\frac{p(e^{tu_1}x_1, ..., e^{tu_n}x_n)}{(e^{tu_1}x_1, ..., e^{tu_n}x_n)^\alpha} = \frac{\sum_\beta c_\beta e^{t \langle u , \beta \rangle}x^\beta}{e^{t\langle u, \alpha \rangle}x^\alpha} = \sum_\beta c_\beta x^{\beta-\alpha} e^{t(\langle u, \beta \rangle - \langle u, \alpha \rangle)}.\]

Taking $t \to \infty$, \[\frac{p(e^{tu_1}x_1, ..., e^{tu_n}x_n)}{(e^{tu_1}x_1, ..., e^{tu_n}x_n)^\alpha} = \sum_\beta c_\beta x^{\beta-\alpha} e^{t(\langle u, \beta \rangle - \langle u, \alpha \rangle) }\to \sum_{\beta \in F} c_\beta x^{\beta-\alpha} = \frac{p_F(x)}{x^\alpha},\]
since $u$ is a normal vector to $F$. 

Hence, for any fixed $x>0,$ \[\inf_{t \ge 0} \frac{p(e^{tu_1}x_1, ..., e^{tu_n}x_n)}{(e^{tu_1}x_1, ..., e^{tu_n}x_n)^\alpha} \le \frac{p_F(x)}{x^\alpha}.\] Thus, \[\text{Cap}_\alpha(p) = \inf_{x>0} \frac{p(x)}{x^\alpha} \le \inf_{x>0} (\inf_{t \ge 0} \frac{p(e^{tu_1}x_1, ..., e^{tu_n}x_n)}{(e^{tu_1}x_1, ..., e^{tu_n}x_n)^\alpha}) \le \inf_{x>0} \frac{p_F(x)}{x^\alpha} =\text{Cap}_\alpha (p_F).\]

\end{proof}

Next, we show that $(s_\lambda)_F$ can be decomposed into the product of Schur polynomials on sets of variables that partition $\{x_1, ..., x_n\}$. 

\begin{lemma}\label{lem:face-schur-decomposition}
    Let $\alpha=(\alpha_1, ..., \alpha_n), \lambda=(\lambda_1, ..., \lambda_n)$ be partitions (allowing zero entries), where the entries of $\lambda$ are not all equal. Suppose $\lambda \trianglerighteq\alpha$ and $\{j \in [n]: \sum_{i=1}^j \lambda_i = \sum_{i=1}^j \alpha_i\} = \{k_1, ..., k_m\}$ where $k_1<...<k_m=n, m>1$. Let $k_0= 0$ and $I_j = \{k_{j-1}+1, ..., k_j\}$ for $j = 1, ..., m$. Then, there exists a face $F$ of $ \text{Newt}(s_\lambda)$ such that $\alpha \in F$ and $(s_\lambda(x))_F = \prod_{j=1}^{m}s_{\lambda_{I_j}}(x_{I_j})$, where $\lambda_{I_j} = (\lambda_{k_{j-1}+1}, ..., \lambda_{k_j})$ and $x_{I_j} = (x_{k_{j-1}+1}, ..., x_{k_j})$. 
\end{lemma}

\begin{proof}
    By \Cref{lem:newton-permutahedron} and the standard facet description of a permutahedron, we can write $\text{Newt}(s_\lambda(x)) = \{\mu : \sum_{i \in T}\mu_i \le \sum_{i=1}^{|T|} \lambda_i, \forall T \subset [n]\}\cap \{\mu: \sum_{i=1}^n \mu_i = \sum_{i=1}^n \lambda_i\}$. Thus, the intersection of $\text{Newt}(s_\lambda(x))$ with the equalities \[\sum_{i=1}^{k_j}\mu_i=\sum_{i=1}^{k_j}\lambda_i,\qquad j=1,\ldots,m,\] is a face $F$ of $\text{Newt}(s_\lambda(x))$, and clearly $\alpha \in F$ by how we defined the indices $k_j$.

    Let $\mathcal{C}$ be the combinatorial class of SSYTs of shape $\lambda$ filled with $1, ..., n$ where each set of consecutive rows with indices in  $I_j$ is filled only with numbers in $I_j$. Let $\mathcal{C}_j$ be the combinatorial class of SSYTs of shape $\lambda_{I_j}$ filled with numbers in $I_j$. Since $s_\lambda(x)$ is the generating function of SSYTs of shape $\lambda$ filled with $1, ..., n$ counted by weight, $(s_\lambda(x))_F$ is the generating function of $\mathcal{C}$ counted by weight. Clearly, for any SSYT in $\mathcal{C}$, we can decompose it by grouping the rows whose indices correspond to each $I_j$ into $m$ SSYTs in $\mathcal{C}_1, ..., \mathcal{C}_{m}$. Conversely, for $m$ SSYTs in $\mathcal{C}_1, ..., \mathcal{C}_{m}$, we can stack them together to get a SSYT in $\mathcal{C}$. The columns of the resulting tableau are strictly increasing because, for each $j<m$, the numbers in $I_{j+1}$ are all greater than those in $I_j$. Thus, we have a weight-preserving bijection between $\mathcal{C}$ and $\mathcal{C}_1 \times ... \times \mathcal{C}_{m}$. It follows that $(s_\lambda(x))_F = \prod_{j=1}^{m} s_{\lambda_{I_j}}(x_{I_j}).$
\end{proof}

With the lemma above, we can decompose $\text{Cap}_\alpha((s_\lambda)_F)$ into a product of the capacities $\text{Cap}_{\alpha_{I_j}}(s_{\lambda_{I_j}}(x_{I_j}))$, where each $\alpha_{I_j}$ is now in the relative interior of the Newton polytope of the corresponding $s_{\lambda_{I_j}}(x_{I_j})$. We can therefore apply the results proved in the previous subsections to each factor. 

\begin{theorem}\label{thm:boundary-capacity-algorithm}
    Let $\alpha, \lambda$ be partitions such that $\lambda \trianglerighteq\alpha$, and the entries of $\lambda$ are not all equal. Let $d=|\lambda|$. Suppose that there exists an integer $j<n$ such that $\alpha_1+...+\alpha_j = \lambda_1+...+\lambda_j$. Let $n \in \mathbb{N}$ such that $n \ge l(\lambda)$. Then, there exists an optimization problem $(P_\alpha)$ with objective function $g_\alpha(y)$ and an algorithm that, given $\epsilon>0$, returns $g^\circ_\alpha$ such that
    \[
    g_\alpha^* \le g_\alpha^\circ\le g_\alpha^* + \epsilon,
    \]
    where $g_\alpha^*$ is the optimal value of ($P_\alpha)$, and
    \[
    \text{Cap}_\alpha(s_\lambda(x_1, ..., x_n)) = e^{g^*_\alpha}.
    \]
    The running time of the algorithm is polynomial in $n, d,$ and $\log(\frac{1}{\epsilon}).$
    
\end{theorem}

\begin{proof}
    First, we can define the face $F$ of $\text{Newt}(s_\lambda)$ as in \Cref{lem:face-schur-decomposition} and decompose $(s_\lambda(x))_F = \prod_{j=1}^{m} s_{\lambda_{I_j}}(x_{I_j})$, where the $I_j$'s partition $[n]$ and are separated at places where the sums of initial parts of $\alpha$ and $\lambda$ are equal. By \Cref{lem:face-capacity}, \[\text{Cap}_\alpha(s_\lambda(x))= \text{Cap}_\alpha((s_\lambda(x))_F) = \prod_{j=1}^{m}\text{Cap}_{\alpha_{I_j}}(s_{\lambda_{I_j}}(x_{I_j})).\]

    Taking the transformation described at the beginning of \Cref{sec:computing-capacity}, we have \[g^*_\alpha = \sum_{j=1}^{m} g^*_{\alpha_{I_j}},\] where $g_\alpha(y) = \ln \sum_{\mu} K_{\lambda, \mu} e^{\langle y, \mu- \alpha\rangle}$, $g_{\alpha_{I_j}}(y_{I_j}) =\ln \sum_{\mu_{I_j}} K_{\lambda_{I_j},\mu_{I_j}} e^{\langle y_{I_j}, \mu_{I_j}-\alpha_{I_j}\rangle}$, and $g^*_\alpha = \inf \{g_\alpha(y): y \in \mathbb{R}^{n}\}, g_{\alpha_{I_j}}^*= \inf \{g_{\alpha_{I_j}}(y_{I_j}): y_{I_j} \in \mathbb{R}^{|I_j|}\}.$

By construction, $\alpha_{I_j} \in \text{relint}(\text{Newt}(s_{\lambda_{I_j}}(x_{I_j})))$, so \Cref{thm:relint-capacity-algorithm} applies. We can run the ellipsoid algorithm on each optimization problem in the decomposition to obtain approximations of the optimal values. Summing these outputs gives the required $g_\alpha^\circ$. Since $m \le n,$ it is enough to approximate each $g^*_{\alpha_{I_j}}$ within $\epsilon/n$. Since $\sum_{j=1}^{m} |\lambda_{I_j}|=|\lambda|= d$ and $\sum_{j=1}^{m} |I_j| = n, $ by \Cref{thm:relint-capacity-algorithm}, the total run time is polynomial in $n, d,$ and $\log(\frac{1}{\epsilon})$.
\end{proof}

\Cref{thm:relint-capacity-algorithm,thm:boundary-capacity-algorithm} show that our ellipsoid method-based algorithms for $(P_\alpha)$ have time complexity polynomial in the same parameters, regardless of whether $\alpha \in \text{relint}(\text{Newt}(s_\lambda)).$

\begin{theorem}\label{thm:capacity-algorithm}
    Let $\alpha, \lambda$ be partitions where $\lambda \trianglerighteq\alpha$, and let $n \ge l(\lambda)$. Let $d = |\lambda|$. There exists an optimization problem $(P_\alpha)$ with objective function $g_\alpha(y)$ and an algorithm that, given $\epsilon>0,$ returns $g_\alpha^\circ$ such that
    \[
    g_\alpha^*\le g_\alpha^\circ \le g^*_{\alpha}+\epsilon,
    \]
    where $g^*_\alpha$ is the optimal value of $(P_\alpha)$, and
    \[
    \text{Cap}_\alpha(s_\lambda(x_1, ...,x_n)) = e^{g^*_\alpha}.
    \]
    The running time of the algorithm is polynomial in $n, d,$ and $\log(\frac{1}{\epsilon})$. 
\end{theorem}

\subsection{A More General Approach}\label{subsec:bit-complexity-approach}

The problem $(P_\alpha)$ can also be viewed as a weighted maximum-entropy dual problem in the framework of \cite{straszak2019maximumentropydistributionsbit}. The geometry of the Newton polytope of $s_\lambda$ satisfies the conditions of their bit-complexity upper bound theorem, which therefore provides an alternative bounded search region containing an $\epsilon$-near-optimal point of $(P_\alpha)$. At the end of this subsection, we compare the resulting bound with the bounding box derived in \cref{subsec:bounding-box}.

\begin{theorem}[Bit Complexity Upper Bound {\cite[Theorem 4.1]{straszak2019maximumentropydistributionsbit}}]\label{thm:sv-bit-complexity}
    Let $\mathcal{F}$ be any finite subset of $\mathbb{Z}^n$, and let $P=\operatorname{conv}(\mathcal{F})$. Then, for every function $p: \mathcal{F}\to \mathbb{R}_{>0}$ and for every number $\epsilon>0$, there exists a number
    \[
    R=O\left(n^{\frac{3}{2}}\cdot M \cdot \left(L_p+\log|\mathcal{F}|+\log\frac{n}{\epsilon}\right)\right)
    \]
    such that 
    \[\forall \theta \in P, \exists y \in B(0, R) \hspace{0.3em} s.t.\hspace{0.3em} h_\theta(y) \le g(\theta)+\epsilon.\]
Here, \begin{itemize}
    \item $h_\theta(y)=\log\sum_{\alpha \in \mathcal{F}}p_\alpha e^{\langle \alpha-\theta, y \rangle}$, where $p_\alpha=p(\alpha)$.
    \item $g(\theta)=\text{inf}_{y \in \mathbb{R}^n}\log\sum_{\alpha \in \mathcal{F}}p_\alpha e^{\langle \alpha-\theta, y \rangle}$.
    \item $M \in \mathbb{N}$, called the \textit{unary facet complexity} of $P$, is the smallest integer such that $P$ has a description of the form $P=\{x\in \mathbb{R}^n: \langle a_i, x\rangle \le b_i, \forall i \in I\}\cap H,$ where $I$ is a finite index set, $a_i \in \mathbb{Z}^n,$ $||a_i||_\infty \le M, $ and $H$ is a linear subspace of $\mathbb{R}^n$ \cite[Definition 4.1]{straszak2019maximumentropydistributionsbit}.
    \item $L_p=\max_{\alpha\in\mathcal{F}}|\log p_\alpha|$.
\end{itemize}
\end{theorem}
     
\begin{theorem}\label{thm:schur-bit-complexity}
    Let $\alpha \in \text{supp}(s_\lambda).$ Let $K(\lambda):=\max \{K_{\lambda, \mu}: \mu \in \operatorname{supp}(s_\lambda) \}$. Then, $\forall \epsilon>0$, there exists a number
    \[
    R=O\left(n^{\frac{3}{2}} \cdot \left(\log K(\lambda)+\log|\text{supp}(s_\lambda)|+\log\frac{n}{\epsilon}\right)\right)
    \]
    such that 
    \[\exists y \in \mathbb{R}^n, ||y|| \le R, g_\alpha(y) \le g_\alpha^*+\epsilon.\]
\end{theorem}
\begin{proof}
    This follows by applying \Cref{thm:sv-bit-complexity}.
    
    Let $u=(d,0,\ldots,0) \in \mathbb{Z}^n$. First, translate $\text{supp}(s_\lambda)$ to $S':=\{\mu-u: \mu \in \text{supp}(s_\lambda)\}$ and define, for $\alpha' \in S'$, \[g'(\alpha'):=\inf_{y\in \mathbb{R}^n} \ln\left(\sum_{\mu' \in S'} K_{\lambda, \mu'+u}e^{\langle \mu'-\alpha', y\rangle}\right).\] Clearly $\forall \alpha\in\text{supp}(s_\lambda), g'(\alpha-u)=g^*_\alpha$. Translating the permutohedron does not change its facet normals, and by \Cref{lem:newton-permutahedron} and the standard permutahedron facet description, it has a facet description with $0/1$ normal vectors. Thus, the unary facet complexity of $\operatorname{conv}(S')$ is 1. The shifted problem is exactly in the form required by \Cref{thm:sv-bit-complexity}: $|S'|=|\text{supp}(s_\lambda)|$ is finite, $S' \subset \mathbb{Z}^n$, and $K_{\lambda, \mu'+u}>0$ for every $\mu' \in S'$. Therefore, plugging in the corresponding parameters gives the result. 
\end{proof}

At first glance, this bounding box for a near-optimal point may look stronger
than our optimizer bounding box in \Cref{subsec:bounding-box}, as it involves
the parameter
$K(\lambda)=\max_{\mu\in \operatorname{supp}(s_\lambda)}K_{\lambda,\mu}$
instead of $\sum_\mu K_{\lambda,\mu}$.
However, this often does not give a better bound in practice, since in most
cases we cannot bound this maximum better than
$
K(\lambda)\le \sum_{\mu}K_{\lambda,\mu}=s_\lambda(1^n).
$
After this substitution, the resulting radius bound is looser, since the
bit-complexity bound also introduces extra factors of $n$. Moreover, the bound
is for an approximate optimizer and depends on the prescribed error. The
advantage is that it applies uniformly to all
$\alpha\in \operatorname{supp}(s_\lambda)$, including boundary points, but for
the main algorithmic analysis we use the Schur-specific bound in the
relative-interior case and the face decomposition for boundary points.

\section{Summary of Results and Numerical Examples}\label{sec:summary-examples}
Combining \Cref{thm:schur-capacity-bound,rem:simple-schur-bound,thm:capacity-algorithm} and using the fact that $\left(\frac{n}{n+|\lambda|}\right)^n \geq e^{-|\lambda|}$, we obtain the main theorem stated in \Cref{thm:main}: the capacity method gives coefficient bounds for nontrivial Kostka numbers, and the relevant capacity can be approximated in polynomial time. Although the preceding results are stated for partitions $\mu$, they extend to arbitrary exponent vectors since Schur polynomials are symmetric.

\begin{theorem}[Main Theorem]
    
\MainTheoremBody

\end{theorem}

As discussed in \Cref{subsec:bit-complexity-approach}, $(P_\mu)$ is an instance of the weighted maximum-entropy dual problems studied in \cite{straszak2019maximumentropydistributionsbit}. Thus, the polynomial dependence of the running time on $\log(1/\epsilon)$ is consistent with their general bit-complexity and algorithmic results. Our analysis verifies the existence of an efficient oracle specifically for Schur polynomials.

Finally, we ran numerical examples to illustrate the efficiency of the capacity method and to compare the resulting bounds with known Kostka numbers and elementary upper bounds. Here the elementary upper bound is taken to be the minimum of $s_\lambda(1^n)$, $f^\lambda$ (the number of standard Young tableaux of shape $\lambda$), and $\binom{|\mu|}{\mu_1,\ldots,\mu_n}$ (the number of words with content $\mu$). The convex optimization was carried out by parametrizing the hyperplane $\sum_i y_i=0$ as $y_n=-\sum_{i=1}^{n-1}y_i$, and then applying BFGS with backtracking line search. The reported residual is the 2-norm of the reduced gradient. We separate face and relative interior cases because the face decomposition in \cref{subsec:boundary-case} can simplify both the capacity and exact computation for face cases. In particular, the exact Kostka number in the large face example was computed using the face decomposition from \Cref{subsec:boundary-case}: direct extraction did not finish within the imposed cutoff. For the larger relative interior examples, direct coefficient extraction also did not finish within the imposed cutoff, so the exact Kostka numbers are not reported. 

\begin{table}[htbp]
\centering
\small
\setlength{\tabcolsep}{3pt}
\renewcommand{\arraystretch}{1.08}
\resizebox{\textwidth}{!}{%
\begin{tabular}{@{}lllrrrr@{}}
\hline
case & $\lambda$ & $\mu$ & $K_{\lambda,\mu}$
& lower bound & $\operatorname{Cap}_{\mu}$
& elem.\ UB \\
\hline
rel.\ int.
& $(2,1,0)$
& $(1,1,1)$
& $2$
& $0.5$
& $8$
& $2$ \\

face
& $(4,3,1,0)$
& $(4,2,1,1)$
& $2$
& $0.7336$
& $11.74$
& $70$ \\

rel.\ int.
& $(8,5,3,1)$
& $(5,5,4,3)$
& $15$
& $1.342$
& $814.2$
& $1260$ \\

face
& $(35,28,21,14,7,0)$
& $(32,31,17,14,11,0)$
& $5$
& $1.54$
& $2175$
& $3.518\times10^{13}$ \\

rel.\ int.
& $(35,28,21,14,7,0)$
& $(18,18,18,17,17,17)$
& {--}
& $1.263\times10^5$
& $3.436\times10^{13}$
& $3.518\times10^{13}$ \\

rel.\ int.
& $(50,40,30,20,10,0)$
& $(25,25,25,25,25,25)$
& {--}
& $2.611\times10^6$
& $4.177\times10^{15}$
& $4.177\times10^{15}$ \\
\hline
\end{tabular}%
}
\caption{Capacity bounds and elementary upper bounds. The lower bound is the sharper product bound from \Cref{thm:schur-capacity-bound}. The exact value in the large face example was obtained by face decomposition; dashes indicate that an exact value was not obtained.}
\label{tab:schur-capacity-bounds}
\end{table}

\begin{table}[htbp]
\centering
\small
\setlength{\tabcolsep}{3pt}
\renewcommand{\arraystretch}{1.08}
\resizebox{\textwidth}{!}{%
\begin{tabular}{@{}lllrrr@{}}
\hline
case & $\lambda$ & $\mu$
& cap.\ time (s)
& exact method/time
& residual \\
\hline
rel.\ int.
& $(2,1,0)$
& $(1,1,1)$
& $2.99\times10^{-4}$
& direct: $1.83\times10^{-4}$ s
& $6.21\times10^{-17}$ \\

face
& $(4,3,1,0)$
& $(4,2,1,1)$
& $0.0115$
& direct: $0.00855$ s
& $9.34\times10^{-12}$ \\

rel.\ int.
& $(8,5,3,1)$
& $(5,5,4,3)$
& $0.958$
& direct: $1.03$ s
& $4.87\times10^{-13}$ \\

face
& $(35,28,21,14,7,0)$
& $(32,31,17,14,11,0)$
& $3.07$
& face decomposition: $1.04\times10^{-3}$ s
& $2.88\times10^{-6}$ \\

rel.\ int.
& $(35,28,21,14,7,0)$
& $(18,18,18,17,17,17)$
& $2.33$
& not obtained
& $6.52\times10^{-5}$ \\

rel.\ int.
& $(50,40,30,20,10,0)$
& $(25,25,25,25,25,25)$
& $0.0398$
& not obtained
& $0$ \\
\hline
\end{tabular}%
}
\caption{Running times for capacity evaluation and exact Kostka-number computation. Here ``direct'' means Schur-to-monomial coefficient extraction, and the residual is the Euclidean norm of the reduced gradient.}
\label{tab:schur-capacity-times}
\end{table}

From the examples, we see that the bounds can behave better when $\mu$ lies on a face of the corresponding Newton polytope. Face decomposition also makes the exact computation here easier. However, this does not mean that the exact computation for face cases can always be done efficiently: the face decomposition may still leave large relative interior subproblems whose exact coefficients are difficult to compute. The computational advantage of the capacity method is that it can still provide efficiently computable bounds in cases where exact coefficient extraction, either directly or after decomposition, is difficult. The numerical experiments also indicate that the capacity can be evaluated efficiently by standard convex optimization methods, with small reduced-gradient residuals. In comparison with other results that motivate non-elementary approaches to approximate Kostka numbers (e.g., see \cite{McAllister2008}, \cite{Narayanan_2026}), which are often asymptotic in nature, our capacity bounds apply directly to fixed pairs $(\lambda,\mu)$ and do not require passing to a dilation. It would be interesting to compare these estimates systematically with existing asymptotic bounds, and to apply the same capacity framework to coefficients of other important Lorentzian or stable polynomial families.

\printbibliography

\end{document}